\documentclass[letterpaper, 10 pt, conference]{ieeeconf}  

\IEEEoverridecommandlockouts                              
\usepackage{epsfig} 
\usepackage{amsmath,amssymb,amsfonts}
\usepackage{color}
\usepackage{booktabs}
\usepackage{cite}
\usepackage{subfigure}
\usepackage{balance}
\usepackage{graphicx}
\usepackage{multirow}
\usepackage{mathtools}
\usepackage{pdfrender}

\usepackage{centernot}
\usepackage{color}
\usepackage{subfig}
\usepackage[acronym]{glossaries}
\usepackage[ruled,vlined]{algorithm2e}
\usepackage{algpseudocode}

\makeglossaries
\newacronym{LP}{LP}{Linear Programming}
\newacronym{psd}{PSD}{positive semidefinite}
\newacronym{SDP}{SDP}{semidefinite program}
\newacronym{FB}{FB}{forward-backward}
\newacronym{pFB}{pFB}{preconditioned forward-backward}
\newglossaryentry{LMI}
{
	name={LMI},
	description={linear matrix inequality},
	first={\glsentrydesc{LMI} (\glsentrytext{LMI})},
	plural={LMIs},
	descriptionplural={linear matrix inequalities},
	firstplural={\glsentrydescplural{LMI} (\glsentryplural{LMI})}
}

\newcommand{\R}{\mathbb{R}}

\newcommand{\N}{\mathbb{N}}

\newcommand{\bbS}{\mathbb{S}}
\newcommand{\mc}[1]{\mathcal{#1}}

\newcommand{\zer}{\mathrm{zer}}

\newcommand{\proj}{\mathrm{proj}}
\newcommand{\diag}{\mathrm{diag}}

\newcommand{\vvec}{\mathrm{vec}}
\newcommand{\mat}{\mathrm{mat}}
\newcommand{\col}{\mathrm{col}}

\newcommand{\inner}[2]{\langle #1, #2 \rangle}

\newcommand{\bs}[1]{\boldsymbol{#1}}
\newcommand{\bsone}{\boldsymbol{1}}

\newcommand{\continuanceref}{}

\newtheorem{lemma}{Lemma}

\newtheorem{example}{Example}

\newtheorem{standing}{Standing Assumption}

\title{\LARGE \bf
A preconditioned Forward-Backward method for partially separable SemiDefinite Programs
}

\author{Filippo Fabiani \and
	 Sergio Grammatico
\thanks{The authors are with the Delft Center for Systems and Control, TU Delft, The Netherlands {\tt\footnotesize (\{f.fabiani, s.grammatico\}@tudelft.nl)}
	This work was partially supported by 3mE/TU Delft under cohesion research project Intelligent Autonomous Vehicles and by the ERC under research project COSMOS (ERC-StG 802348).}%
}

\begin{document}

\maketitle
\thispagestyle{empty}
\pagestyle{empty}

\begin{abstract}
We present semi-decentralized and distributed algorithms, designed via a preconditioned forward-backward operator splitting, for solving large-scale, decomposable \glspl{SDP}. 
We exploit a chordal aggregate sparsity pattern assumption on the original \gls{SDP} to obtain a set of mutually coupled \glspl{SDP} defined on \gls{psd} cones of reduced dimensions. 
We show that the proposed algorithms converge to a solution of the original \gls{SDP} via iterations of reasonable computational cost. Finally, we compare the performances of the two proposed algorithms with respect to others available in the literature.
\end{abstract}

\section{Introduction}
Semidefinite programs (SDPs), i.e., convex optimization problems over the cone of \gls{psd} matrices, find applications in many engineering areas, since semidefinite constraints and linear matrix inequalities characterize problems related to sensor networks \cite{biswas2006semidefinite,simonetto2014distributed}, smart grids \cite{dall2013distributed}, robustness analysis of interconnected, uncertain systems \cite{andersen2014robust}, machine learning and operations research \cite{boyd1994linear}. Additionally, \gls{SDP} offer a convex relaxation framework that is applicable to several nonconvex optimization problems, e.g., optimal power flow \cite{madani2014convex,magnusson2015distributed}.

Despite these appealing features, the computational complexity of \glspl{SDP} represents an important issue, since they scale badly with both the dimension of the semidefinite constraints and the number of optimization variables. This dramatically limits, or even makes vain \cite{andersen2014robust, kim2009exploiting}, the possibility of solving large-scale problems in a centralized fashion. Therefore, since certain applications impose intrinsic privacy limitations that rule out centralized computation, it is fundamental to develop scalable and parallelizable algorithms.

Existing works that propose distributed approaches rely on a specific sparsity structure of the \gls{SDP}. Specifically, when the graph associated with the aggregate sparsity pattern of the data matrices enjoys the chordality property, the original \gls{SDP} can be decomposed into a set of mutually coupled \glspl{SDP} defined on cones of smaller dimensions. 
However, the available distributed, first-order algorithms are either ad-hoc for the specific problem considered, e.g. \cite{dall2013distributed,weng2013distributed,simonetto2014distributed,biswas2006semidefinite}, or are variants of the alternating direction method of multipliers (ADMM) \cite{kalbat2015fast,madani2017low}, a special instance of the Douglas-Rachford operator splitting \cite[\S 25.2]{bauschke2011convex}. Along the same line, in \cite{sun2014decomposition} an operator-splitting method has been proposed, whose convergence requires to solve a quadratic program at every iteration. The authors in \cite{nayakkankuppam2007solving} introduced a proximal bundle method encompassing a message-passing strategy, which in turn has been recently exploited in \cite{pakazad2017distributed} together with an interior-point method, to overcome the typical issues of second-order methods shown, e.g., in \cite{fukuda2001exploiting,gros2014newton}.

Differently from the aforementioned literature, in this paper we provide a unifying operator-theoretic perspective to the design of scalable and distributed algorithms for decomposable, large-scale \glspl{SDP}, paving the way for a systematic development of efficient methods that rely on monotone operator theory as a general purpose tool. Specifically, we discuss a semi-decentralized and a fully distributed version of a \gls{pFB} method, recently proposed in \cite{yi2019operator} to solve generalized Nash equilibrium problems. Unlike the standard \gls{FB} method, a preconditioning matrix allows to circumvent intrinsic technical issues, e.g., the impossibility to explicitly compute the resolvent of an operator, by determining the step sizes and the exploration direction of the forward iteration. We exploit the chordal structure of the original problem to obtain a multi-agent \gls{SDP} with coupling constraints. In this context, the proposed algorithms guarantee global convergence to an optimal solution by requiring local evaluations of projections on \gls{psd} cones of reduced dimensions, which is the major computational burden.

After briefly reviewing some notions of graph theory (\S II), under the assumption of chordal sparsity pattern we describe how to decompose an \gls{SDP} into mutually coupled \glspl{SDP} (\S III). Finally, in \S IV, we apply an operator-splitting approach for the design of the \gls{pFB}-based algorithms, and in  \S V, we compare the two methods with available algorithms in the literature on illustrative examples.

\subsubsection*{Basic notation} 
$\R^n$ denotes the set of $n \times 1$ real vectors, $\bbS^{n}$ is the space of $n \times n$ symmetric matrices and $\bbS_{\succcurlyeq 0}^{n} \coloneqq \left\{X \in \bbS^n \mid X \succcurlyeq 0\right\}$ is the cone of \gls{psd} matrices. $\bs{0}$ ($\bsone$) denotes a matrix/vector with all elements equal to $0$ ($1$). 
For vectors $v_1,\dots,v_N\in\mathbb{R}^n$ and $\mc I=\{1,\dots,N \}$, we denote $\boldsymbol{v}:= (v_1 ^\top,\dots ,v_N^\top )^\top = \mathrm{col}(\{v_i\}_{i\in\mc I})$.
The symbol $\inner{\cdot}{\cdot}$ denotes the inner product in the appropriate space, i.e., $\inner{x}{y} = x^\top y$ for $x$, $y \in \R^n$ and $\inner{X}{Y} = \textrm{trace}\left(XY\right)$ for $X$, $Y \in \bbS^n$. Given a matrix $A \in \R^{m \times n}$, its $(i,j)$ entry is denoted by $[A]_{i,j}$, while $A^\top$ denotes its transpose. $A \otimes B$ represents the Kronecker product between the matrices $A$ and $B$. Given a discrete set $\mc{X} \subseteq \R^{n}$, $|\mc{X}|$ denotes its cardinality. The operator $\vvec : \bbS^n \to \R^{n^2}$ maps a matrix to a vector that stacks its columns; $\mat(\cdot) \coloneqq \vvec^{-1}(\cdot)$ performs the inverse.
\subsubsection*{Operator-theoretic definitions} 
 $\textrm{Id}(\cdot)$ denotes the identity operator. The mapping $\iota_{\mc{X}} : \R^n \to \{0, \infty\}$ represents the indicator function for a set $\mc{X} \subseteq \R^n$, i.e., $\iota_{\mc{X}}(x) = 0$ if $x \in \mc{X}$, $\infty$ otherwise. The mapping $\proj_{\mc{X}}(\cdot)$ denotes the projection onto the set $\mc{X}$ w.r.t. the Euclidean norm when the argument is a vector, to the Frobenius norm when it is a matrix. The set-valued mapping $\textrm{N}_{\mc{X}} : \R^n \rightrightarrows \R^n$ denotes the normal cone operator for the set $\mc{X}$, i.e., $\textrm{N}_{\mc{X}}(x) = \varnothing$ if $x \notin \mc{X}$, $\{v \in \R^n \mid \sup_{z \in \mc{X}} v^\top (z - x) \leq 0\}$ otherwise. $\zer(\mc{F}) \coloneqq \{x \in \R^n \mid \bs{0} \in \mc{F}(x)\}$ denotes the set of zeros of a set-valued mapping $\mc{F}:\R^n \rightrightarrows \R^n$.

\section{Mathematical problem setup}
We first recall some notions of graph theory and chordality \cite{blair1993introduction,mesbahi2010graph}, necessary to introduce the partially separable structure addressed in the remainder. Then, we review a key result to decompose a convex cone, i.e., the cone of \gls{psd} matrices, into a set of smaller but coupled convex cones. 

\subsection{Graph theoretic concepts}
Let $\mc{G}\coloneqq(\mc{V},\mc{E})$ be a graph defined by a set of vertices $\mc{V} \coloneqq \left\{1,2,\ldots,n\right\}$ and a set of edges $\mc{E} \subseteq \mc{V} \times \mc{V}$. The graph $\mc{G}$ is undirected when $(i,j) \in \mc{E}$ if and only if $(j,i) \in \mc{E}$, while it is complete if every two vertices share exactly one edge. The set of neighboring vertices of $i \in \mc{V}$ is denoted by $\mc{N}_i \coloneqq \{j \in \mc{V} \mid (j,i) \in \mc{E} \}$. A path of $\mc{G}$ is a sequence of distinct vertices such that any consecutive vertices in the sequence correspond to an edge of graph $\mc{G}$. Let $W \in \R^{n \times n}$ be the adjacency matrix of $\mc{G}$ with $[W]_{i,j} \eqqcolon w_{i,j} > 0$ if $j \in \mc{N}_i$, $w_{i,j} = 0$ otherwise, and let us assume $W \in \bbS^n$. The weighted degree matrix is $\Delta \coloneqq \diag(\{d_i\}_{i \in \mc{V}})$, $d_i \coloneqq \sum_{j \in \mc{V}} w_{i,j}$, and the weighted Laplacian of $\mc{G}$ is $L \coloneqq \Delta - W$. 

A clique is defined as the complete subgraph induced by the set of vertices $\mc{C} \subseteq \mc{V}$, i.e., such that $(i,j) \in \mc{E}$ for any distinct vertices $i, j \in \mc{C}$, and it is not a subset of any other clique. 
Given a path $\left\{v_1, v_2, \ldots, v_k\right\}$, a chord is an edge $\left(v_i, v_j\right)$ with $|j - i| > 1$, while a cycle is a path with $(v_1, v_k) \in \mc{E}$ and $(v_i, v_{i+1}) \in \mc{E}$ for $i = 1, 2,\ldots, k-1$. An undirected graph $\mc{G}$ is chordal if every cycle with $k > 3$ has at least one chord. Every nonchordal graph can be extended to a chordal graph by adding a (minimum) number of edges to the original graph (minimum fill-in problem, NP-complete \cite{blair1993introduction}).
Typical examples of chordal graphs are the complete graphs, trees and forests, and undirected graphs with no cycles of length greater than three.

\subsection{Separable cones and chordal decomposition}
Given a sparse matrix $X \in \bbS^n$, the associated sparsity pattern can be represented as an undirected graph $\mc{G}\coloneqq(\mc{V},\mc{E})$ where each clique $\mc{C}_i$ defines a maximally dense principal submatrix of $X$, i.e., $X_i = \mc{P}_{i}(X) \coloneqq P_{\mc{C}_{i}} X P_{\mc{C}_{i}}^\top \in \bbS^{|\mc{C}_i|}$, by means of an entry selection matrix $P_{\mc{C}_{i}} \in \R^{|\mc{C}_{i}| \times n}$ defined as
$$
\left[P_{\mc{C}_{i}}\right]_{h,k} \coloneqq \left\{
\begin{aligned}
&1 \; \text{ if } \mc{C}_{i}(h) = k,\\
&0 \; \text{ otherwise},
\end{aligned}
\right.
$$
where $\mc{C}_{i}(h)$ is the $h$-th element in $C_{i}$, sorted in natural ordering. Conversely, $Y_i = \mc{P}^{-1}_i(X_i) \coloneqq P^\top_{\mc{C}_{i}} X_i P_{\mc{C}_{i}} \in \bbS^{n}$ returns a matrix in the original space. Note that, given cliques $\mc{C}_i$ and $\mc{C}_{j}$ associated to the same graph, $X_i$ and $X_{j}$ have overlapping elements if $\mc{C}_i \cap \mc{C}_{j} \neq \varnothing$.

By considering $\bar{\mc{E}} \coloneqq \mc{E} \cap \left\{(i,i) \mid i \in \mc{V} \right\}$, we define the space of symmetric matrices with sparsity pattern $\mc{G}$ as
$$
\bbS^n(\mc{E}) \coloneqq \left\{ X \in \bbS^n \mid \left[X\right]_{i,j} = \left[X\right]_{j,i} = 0 \text{ if } (i,j) \notin \bar{\mc{E}}\right\}.
$$
Next, we recall a well-known result that links a matrix $X$ belonging to the cone generated by the projection of the \gls{psd} cone onto the space of matrices with sparsity pattern $\mc{G}$, i.e., $X \in \proj_{\bbS^n(\mc{E})} \left(\bbS_{\succcurlyeq 0}^{n}\right)$, and its maximal principal submatrices.
\begin{lemma}\textup{\cite[Th.~7]{grone1984positive}}\label{th:grone}
	Let $\mc{G}$ be a chordal graph and let $\left\{\mc{C}_1, \ldots, \mc{C}_N\right\}$ be the set of its cliques. Then, $X \in \proj_{\bbS^n(\mc{E})} \left(\bbS_{\succcurlyeq 0}^{n}\right)$ if and only if
	$$\mc{P}_{{i}}(X) \in \bbS^{|\mc{C}_i|}_{\succcurlyeq 0}, \; \text{ for all }  i \in \left\{1,\ldots,N\right\}.$$
	\hfill$\square$
\end{lemma}
Informally speaking, a matrix $X \in \bbS^n(\mc{E})$ is positive semidefinite if and only if its maximal dense principal submatrices $X_i \in \bbS^{|\mc{C}_i|}$, $i = 1, 2, \ldots,N$, are positive semidefinite, i.e., the \gls{psd} property of $X$ can be checked on $N$ matrices of smaller dimensions.
To recap the introduced concepts, let us show an example on a simple graph.
\begin{example}
	\begin{figure}
	\centering
	\subfigure[]{%
		\includegraphics[width=0.49\columnwidth]{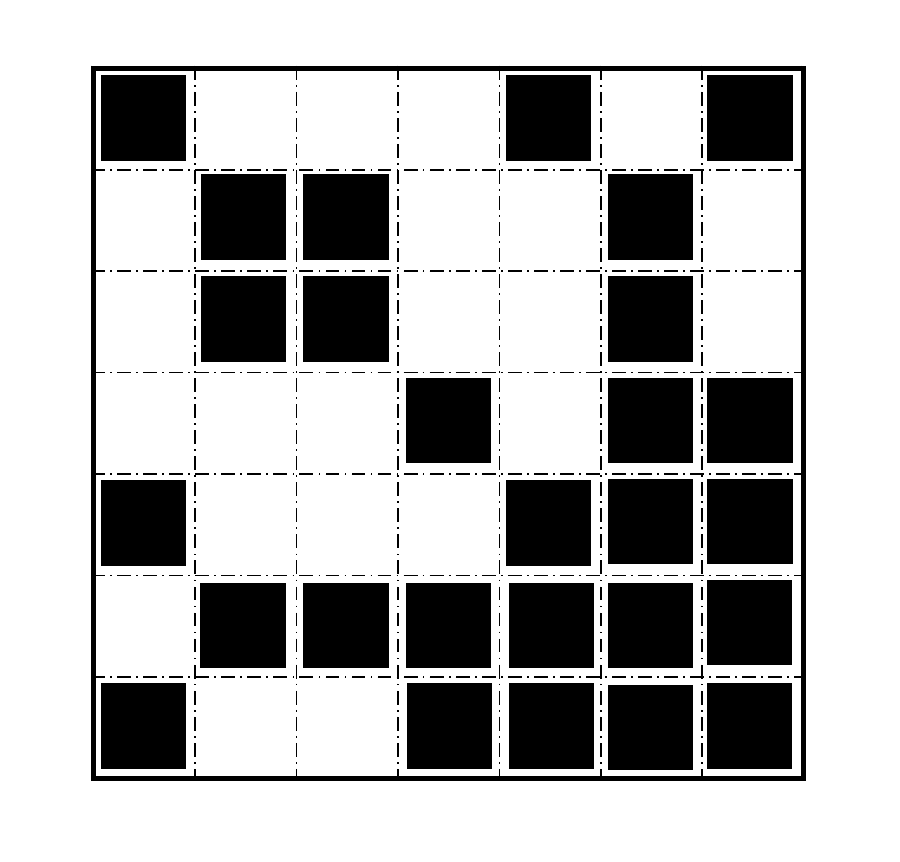}}
	\hfill
	\subfigure[]{%
		\includegraphics[width=0.49\columnwidth]{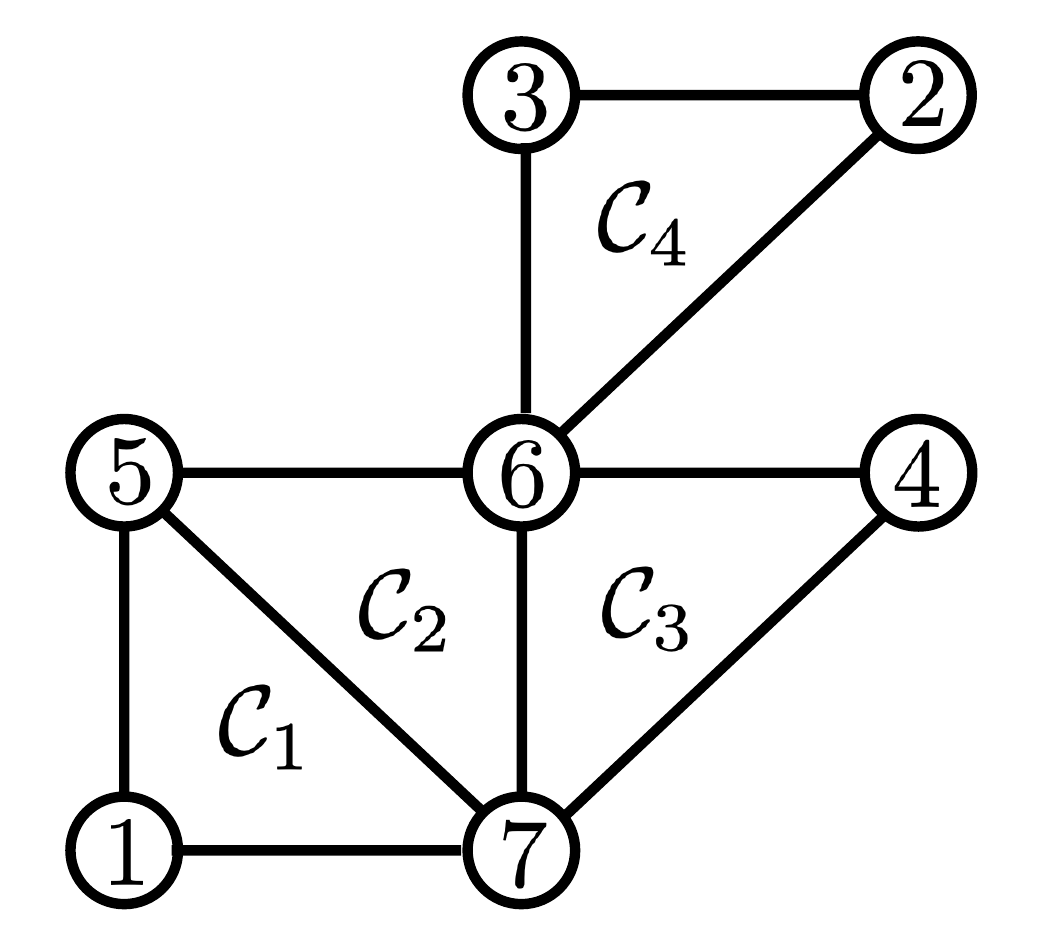}}
		\caption{(a) Sparsity pattern of a given matrix $X \in \bbS{^7}$. (b) Graph associated with the sparsity pattern of $X$.}
		\label{fig:chordal_ex}
	\end{figure}

	The sparsity pattern of a given matrix $X \in \bbS^7$ in Fig.~\ref{fig:chordal_ex}(a) can be summarized in the graph $\mc{G}$ with seven vertices depicted in Fig.~\ref{fig:chordal_ex}(b). The graph $\mc{G}$ happens to be chordal, since every cycle of length strictly greater than three has a chord, e.g., $(6,7)$ is a chord for the cycle $\{5,7,4,6,5\}$. Moreover, it is characterized by four cliques, namely $\mc{C}_1 = \{1,5,7\}$, $\mc{C}_2 = \{5,6,7\}$, $\mc{C}_3 = \{4,6,7\}$ and $\mc{C}_4 = \{2,3,6\}$. For instance, let us consider $\mc{C}_1$. Then,
	$$
	\begin{aligned}
	P_{\mc{C}_1} &= \left[\begin{array}{ccccccc}
	1 & 0 & 0 & 0 & 0 & 0 & 0\\
	0 & 0 & 0 & 0 & 1 & 0 & 0\\
	0 & 0 & 0 & 0 & 0 & 0 & 1\\
	\end{array}\right],\\
	X_1 &= \left[\begin{matrix}
	[X]_{1,1} & [X]_{1,5} & [X]_{1,7}\\
	[X]_{5,1} & [X]_{5,5} & [X]_{5,7}\\
	[X]_{7,1} & [X]_{7,5} & [X]_{7,7}
	\end{matrix}\right],
	\end{aligned}
	$$
	where $X_1$ is the maximally dense principal submatrix of $X$ related with $\mc{C}_1$. It follows from Lemma~\ref{th:grone} that $X$ is \gls{psd} with sparsity pattern $\mc{G}$ if and only if $X_i \in \bbS^{|\mc{C}_i|}_{\succcurlyeq 0}$, $i \in \{1,2,3,4\}$.
	Finally, note that $\mc{C}_1 \cap \mc{C}_2 = \{5,7\}$, and therefore matrices $X_1$ and $X_2$ shall be coupled via the elements $[X]_{5,7} = [X]_{7,5}$, $[X]_{5,5}$ and $[X]_{7,7}$. The same considerations apply to  $\mc{C}_2 \cap \mc{C}_3 = \{6,7\}$ and $\mc{C}_2 \cap \mc{C}_3 \cap \mc{C}_4 = \{6\}$.
	\hfill$\square$
\end{example}

\section{Decomposition in sparse \glspl{SDP}}
We consider an \gls{SDP} optimization problem in standard primal form \cite{vandenberghe2015chordal}:
\begin{equation}\label{eq:SDP}
\left\{
\begin{aligned}
&\underset{X}{\textrm{min}} & & \inner{C}{X}\\
&\textrm{ s.t. } & & \inner{A_{k}}{X} = b_k, \; \forall k \in \mc{K},\\
&&& X \in \bbS_{\succcurlyeq 0}^{n},
\end{aligned}	
\right.
\end{equation}
where $C \in \bbS^n$, $A_k \in \bbS^n$ and $b_k \in \R$ for all $k \in \mc{K} \coloneqq \left\{1,2,\ldots,m\right\}$, whose optimal solutions belong to the set
\begin{equation}\label{eq:SDP_argmin}
\bs{X}^\star \coloneqq 
\underset{X \in \bbS_{\succcurlyeq 0}^{n}}{\textrm{argmin}} \; \inner{C}{X} \textrm{ s.t. } \inner{A_{k}}{X} = b_k, \; \forall k \in \mc{K},
\end{equation}
which is assume to be nonempty. Moreover, we make the following assumption that characterize the data in \eqref{eq:SDP}.
\begin{standing}
	The aggregate sparsity pattern of \eqref{eq:SDP} is identified by the graph $\mc{G} \coloneqq \left(\mc{V}, \mc{E}\right)$, i.e., $C \in \bbS^n(\mc{E})$ and $A_k \in \bbS^n(\mc{E})$, for all $k \in \mc{K}$.
	\hfill$\square$
\end{standing}
Such a sparsity pattern also reflects on the standard dual form of \eqref{eq:SDP} with dual variables $y \in \R^m$ and $Z \in \bbS^n$, i.e.,
\begin{equation}\label{eq:SDP_dual}
\left\{
\begin{aligned}
&\underset{y,Z}{\textrm{max}} & & \inner{b}{y}\\
&\textrm{ s.t. } & & Z + \sum_{k \in \mc{K}} A_k y_k = C,\\
&&& Z \in \bbS_{\succcurlyeq 0}^{n},
\end{aligned}	
\right.
\end{equation}
where $b \coloneqq \col\left(\{b_k\}_{k \in \mc{K}}\right) \in \R^m$. 
Clearly, in view of the equality constraint in \eqref{eq:SDP_dual}, the slack matrix $Z$ is necessarily sparse at any dual feasible point, with sparsity pattern $\mc{G}$. On the other hand, the primal variable $X$ is typically dense, but one can observe that both the cost function and the equality constraints in \eqref{eq:SDP} depend on the entries $\left[X\right]_{i,j}$ such that $(i,j) \in \bar{\mc{E}}$. The remaining entries are arbitrary and guarantee the \gls{psd} property of $X$ (matrix completion problem). Thus, according to the domain-space decomposition in \cite{fukuda2001exploiting}, the conic constraint in \eqref{eq:SDP} can be equivalently rewritten as $X \in \proj_{\bbS^n(\mc{E})} \left(\bbS_{\succcurlyeq 0}^{n}\right)$. We then assume the following.
\begin{standing}\label{standing:chordal}
	The graph $\mc{G}$ is chordal, undirected and connected.
	\hfill$\square$
\end{standing}

Thus, in view of Lemma~\ref{th:grone}, we follow the so called clique-tree conversion introduced in \cite{fukuda2001exploiting,kim2011exploiting} to decompose the original \gls{SDP} in \eqref{eq:SDP}. Specifically, let $\left\{\mc{C}_i\right\}_{i \in \mc{I}}$, $\mc{I} \coloneqq \left\{1,2,\ldots,N\right\}$, be the set of cliques of $\mc{G}$, and let $\mc{N}_i$ be the set containing the indices of all the cliques such that $\mc{C}_i \cap \mc{C}_j \neq \varnothing$, $j \in \mc{I}$. The \gls{SDP} in \eqref{eq:SDP} translates into
\begin{equation}\label{eq:SDP_decomposed}
\left\{
\begin{aligned}
&\underset{X_1, \ldots, X_N}{\textrm{min}} & &\hspace{-.3cm} \sum_{i \in \mc{I}} \inner{C_i}{X_i}\\
&\hspace{.3cm}\textrm{ s.t. } & &\hspace{-.3cm} \sum_{i \in \mc{I}} \inner{A_{k,i}}{X_i} = b_k, \; \forall k \in \mc{K},\\
&&&\hspace{-.3cm} \mc{P}_{i,j} \left( \mc{P}^{-1}_{i}(X_i) - \mc{P}^{-1}_{j}(X_j) \right) = 0, \; \forall (i,j) \in \mc{I} \times \mc{N}_i,\\
&&&\hspace{-.3cm} X_i \in \bbS_{\succcurlyeq 0}^{|\mc{C}_i|}, \; \forall i \in \mc{I},
\end{aligned}	
\right.
\end{equation}
where $\mc{P}_{i,j}(X) = \mc{P}_{j,i}(X) \coloneqq P_{\mc{C}_{i} \cap \mc{C}_{j}} X P_{\mc{C}_{i} \cap \mc{C}_{j}}^\top \in \bbS^{|\mc{C}_i \cap \mc{C}_{j}|}$, while the coefficient matrices $C_i$ and $A_{k,i}$ are chosen so that $\inner{C}{Q} = \sum_{i \in \mc{I}} \inner{C_i}{Q_{i}}$ and $\inner{A_k}{Q} = \sum_{i \in \mc{I}} \inner{A_{k,i}}{Z_{i}}$ for all $Q \in \bbS^n(\mc{E})$. Remarkably, their choice is not unique. Some examples on how to chose such matrices are in \cite{sun2014decomposition}, while for the remainder we postulate the following assumption.
\begin{standing}\label{ass:independence}
	For all $i \in \mc{I}$, $A_{1,i}, \ldots, A_{m,i}$ are linearly independent matrices. 
	\hfill$\square$
\end{standing}

The \gls{SDP} in \eqref{eq:SDP_decomposed} is fully decomposed in terms of cost function and it is defined on \gls{psd} conic constraints of smaller dimensions. However, it shows couplings among the dense submatrices in the two sets of equality constraints. Specifically, while the first set is directly inherited from \eqref{eq:SDP}, the second set refers to the consistency constraints, which ensure that the entries of the dense submatrices $X_i$ and $X_j$ agree when share the same entries of the original $X$ in \eqref{eq:SDP}. Note that, all these constraints are affine in the variable $X_i$, since they involve matrix multiplications only. Thus, the decomposed structure in \eqref{eq:SDP_decomposed}, which is quite general for a sparse \gls{SDP} under Standing Assumption~\ref{standing:chordal}, enables and motivates us to develop semi-decentralized/distributed algorithms that allow to compute a solution $X^\ast \in \bs{X}^\star$ by treating each clique $i \in \mc{I}$ as a single decision agent that controls its local variable, $X_i \in \bbS^{n_i}$. 

\subsection{Vector notation for decomposed \glspl{SDP}}
For our purposes, i.e., to formalize our multi-agent optimization problem and, consequently, to develop solution algorithms,  it might be convenient to adopt a vector notation for the matrix variables in \eqref{eq:SDP_decomposed}. Specifically, let us now define $n_i \coloneqq |\mc{C}_i|$, for all $i \in \mc{I}$, and $n_{i,j} \coloneqq |\mc{C}_i \cap \mc{C}_j|$ for all $i \in \mc{I}$ and all $j \in \mc{N}_i$. The \gls{SDP} in \eqref{eq:SDP_decomposed} can be rewritten as multi-agent \gls{SDP} and recast in vector form as follows:
\begin{equation}\label{eq:SDP_single_prob_vecform}
\left\{
\begin{aligned}
&\underset{x_1, \ldots, x_N}{\textrm{min}} & & \sum_{i \in \mc{I}} \inner{c_i}{x_i}\\
&\hspace{.25cm}\textrm{ s.t. } & & A_i x_i + \sum_{j \in \mc{I} \setminus \{i\}} A_j x_j = b,\\
&&& D_{i,j} x_i +  D_{j,i} x_j = 0, \; \forall (i,j) \in \mc{I} \times \mc{N}_i,\\
&&& x_i \in \mc{S}_i, \; \forall i \in \mc{I},
\end{aligned}	
\right.
\end{equation}
where, for all $i \in \mc{I}$, $x_i \coloneqq \vvec(X_i) \in \R^{n_i^2}$, $c_i \in \vvec(C_i) \in \R^{n_i^2}$, $A_i \coloneqq [\vvec(A_{1,i}) \; \vvec(A_{2,i}) \; \ldots \; \vvec(A_{m,i}) ]^\top \in \R^{m \times n_i^2}$, which is of full-row rank in view of Assumption~\ref{ass:independence}. Moreover, $D_{i,j} \coloneqq P_{i,j}  P_i \in \R^{n_{i,j} \times n_i^2}$ and $D_{j,i} \coloneqq - P_{i,j}  P_j \in \R^{n_{i,j} \times n_j^2}$, where $P_i \coloneqq \left(P_{\mc{C}_i} \otimes P_{\mc{C}_i}\right)^\top \in \R^{n \times n_i^2}$, $n \coloneqq \sum_{i \in \mc{I}} n_i^2$, for all $i \in \mc{I}$ and $P_{i,j} = P_{j,i} \coloneqq P_{\mc{C}_{j} \cap \mc{C}_{i}} \otimes P_{\mc{C}_{j} \cap \mc{C}_{i}} \in \R^{n_{i,j} \times n}$. The set $\mc{S}_i$ is defined as
$
\mc{S}_i \coloneqq \{x_i \in \R^{n_i^2} \mid \mat(x_i) \in \bbS^{n_i}_{\succcurlyeq 0} \}.
$
It follows directly that $x_i \in \mc{S}_i$ if and only if $X_i \in \bbS^{n_i}_{\succcurlyeq 0}$. 

\section{Preconditioned Forward-Backward Operator Splitting}
In this section, we look at the multi-agent \gls{SDP} in \eqref{eq:SDP_decomposed} or, equivalently, in \eqref{eq:SDP_single_prob_vecform}, under an operator splitting perspective. By leveraging on this view, we then provide two iterative algorithms based on a preconditioned variant of the \gls{FB} operator splitting \cite[\S 25.3]{bauschke2011convex}.

Let us introduce $\bs{x} \coloneqq \col\left( \left\{x_i\right\}_{i \in \mc{I}} \right) \in \R^{n}$. The Lagrangian function associated to \eqref{eq:SDP_single_prob_vecform} reads as follows
$$
\begin{aligned}
\mc{L} (\bs{x}, \nu, \lambda) \coloneqq &\textstyle\sum_{i \in \mc{I}} \inner{c_i}{x_i} + \iota_{\mc{S}}(\bs{x}) + \nu^\top\left(A \bs{x} - b\right)\\ 
+ &\textstyle\sum_{i \in \mc{I}}\textstyle\sum_{j \in \mc{N}_i} \lambda^\top_{i,j} \left(D_{i,j} x_i + D_{j,i} x_j\right),
\end{aligned}
$$
where $\mc{S} \coloneqq \prod_{i \in \mc{I}} \mc{S}_i$, $A \coloneqq [A_1 \, \ldots \, A_N]$, $\nu \in \R^m$ and $\lambda_{i,j} \in \R^{n_{i,j}}$ are the Lagrangian multipliers associated with the coupling constraints. Moreover, let us define $\lambda_{i} \coloneqq \col( \left\{ \lambda_{i,j} \right\}_{j \in \mc{N}_i} ) \in \R^{m_i}$, $m_i \coloneqq \sum_{j \in \mc{N}_i} n_{i,j}$, and $\lambda \coloneqq \col( \left\{ \lambda_{i} \right\}_{i \in \mc{I}} ) \in \R^p$, $p \coloneqq \sum_{i \in \mc{I}} m_i$. 
Then, for each agent $i \in \mc{I}$, we define the consistency constraint matrix as $D_i \in \R^{p \times n^2_i}$ by stacking and rearranging all $D_{i,j}$ and $-D_{i,j}$, together with some (possibly rectangular) matrices of zeros with suitable dimensions.  Thus, to compute a solution to the \gls{SDP} in \eqref{eq:SDP_single_prob_vecform}, the KKT stationarity conditions require that
\begin{equation}\label{eq:KKT}
\left\{
\begin{aligned}
&0 \in \nabla_{\bs{x}} \textstyle \sum_{i \in \mc{I}}\inner{c_i}{x_i} + \textsc{N}_{\mc{S}}(\bs{x}) + A^\top \nu + D^\top \lambda,\\
&0 = - A \bs{x} + b,\\
&0 = - D \bs{x}.
\end{aligned}
\right.
\end{equation}
where $D \coloneqq [D_1 \, \ldots \, D_N]$. To cast the inclusion problem \eqref{eq:KKT} in a compact form, we first introduce the gradient mapping, which is constant in this case, as $F \coloneqq \col(\left\{\nabla_{x_i} \inner{c_i}{x_i} \right\}_{i \in \mc{I}} )$. Then, we define the set-valued mapping $T : \R^{n} \times \R^m \times \R^p  \rightrightarrows \R^{n} \times \R^m \times \R^p$, as follows:
\begin{equation}\label{eq:mapping}
	T : \left[\begin{array}{c}
	\bs{x}\\
	\nu\\
	\lambda
	\end{array}\right] \mapsto \left[\begin{array}{c}
	\textsc{N}_{\mc{S}}(\bs{x}) + F + A^\top \nu + D^\top \lambda\\
	- A \bs{x} + b\\
	- D \bs{x}
	\end{array}\right].
\end{equation}

Thus, solving \eqref{eq:KKT} is equivalent to finding $\left(\bs{x}, \nu, \lambda \right)$ such that $\bs{0} \in T\left(\bs{x}, \nu, \lambda \right)$, i.e., finding a zero of $T$ in \eqref{eq:mapping}, which can be written as the sum of the two following operators:
\begin{align}
\mc{A} : &\left[\begin{array}{c}
\bs{x}\\
\nu\\
\lambda
\end{array}\right] \mapsto  \left[\begin{array}{c}
F\\
b\\
0
\end{array}\right],\label{eq:A}\\
\mc{B} : &\left[\begin{array}{c}
\bs{x}\\
\nu\\
\lambda
\end{array}\right] \mapsto \left[\begin{array}{c}
\textsc{N}_{\mc{S}}(\bs{x})\\
0\\
0
\end{array}\right] + \left[\begin{array}{ccc}
\phantom{-}0 & A^\top & D^\top\\
-A & 0 & 0\\
- D & 0 & 0
\end{array}
\right] \left[\begin{array}{c}
\bs{x}\\
\nu\\
\lambda
\end{array}\right].\label{eq:B}
\end{align}

We note that the set of zeros of $T \coloneqq \mc{A} + \mc{B}$ corresponds to the set of zeros of the operator $\Phi^{-1} \left(\mc{A} + \mc{B}\right)$, where $\Phi \succ 0$ is a ``preconditioning'' matrix to be designed freely. In simple words, the \gls{FB} operator splitting states that the zeros of $\mc{A} + \mc{B}$, or, equivalently, $\Phi^{-1}(\mc{A} + \mc{B})$, are the fixed points of the operator $\left(\textrm{Id} + \Phi^{-1} \mc{B}\right)^{-1} \circ \left(\textrm{Id} - \Phi^{-1} \mc{A}\right)$. Consequently, by introducing $\bs{\omega} \coloneqq \col(\bs{x}, \nu, \lambda)$, the \gls{FB} algorithm corresponds to the fixed-point iteration \cite[(1.67)]{bauschke2011convex}
\begin{equation}\label{eq:PB}
	\bs{\omega}^{k+1} = \underbrace{\left(\textrm{Id} + \Phi^{-1} \mc{B}\right)^{-1}}_{\text{backward}} \circ \underbrace{\left(\textrm{Id} - \Phi^{-1} \mc{A}\right)}_{\text{forward}} \left(\bs{\omega}^{k}\right).
\end{equation}

\subsection{Semi-decentralized forward-backward algorithm}
From \eqref{eq:PB}, we obtain the following chain of inclusions:
$$
\begin{aligned}
\left( \textrm{Id} - \Phi^{-1} \mc{A} \right)& \left( \bs{\omega}^{k} \right) \in \left( \textrm{Id} + \Phi^{-1} \mc{B} \right) \left( \bs{\omega}^{k+1} \right) \iff\\
& - \mc{A} \left( \bs{\omega}^{k} \right) \in \Phi \left( \bs{\omega}^{k+1} - \bs{\omega}^{k}\right) + \mc{B} \left( \bs{\omega}^{k+1} \right),
\end{aligned}
$$
and by explicitly including the operators in \eqref{eq:A}--\eqref{eq:B} we have
\begin{align}
&-\left[\begin{array}{c}
F\\
b\\
0
\end{array}\right] \in \left[\begin{array}{ccc}
\Phi_{11} & \Phi_{12} & \Phi_{13}\\
\Phi_{21} & \Phi_{22} & \Phi_{23}\\
\Phi_{31} & \Phi_{32} & \Phi_{33}
\end{array}
\right] \left[\begin{array}{c}
\bs{x}^{k+1} - \bs{x}^{k}\\
\nu^{k+1} - \nu^{k}\\
\lambda^{k+1} - \lambda^{k}
\end{array}\right]\nonumber\\
&+ \left[\begin{array}{c}
\textsc{N}_{\mc{S}}(\bs{x}^{k+1})\\
0\\
0
\end{array}\right] + \left[\begin{array}{ccc}
\phantom{-}0 & A^\top & D^\top\\
-A & 0 & 0\\
- D & 0 & 0
\end{array}
\right] \left[\begin{array}{c}
\bs{x}^{k+1}\\
\nu^{k+1}\\
\lambda^{k+1}
\end{array}\right].\label{eq:low_term}
\end{align}

Thus, by following the guidelines provided in \cite[\S IV]{belgioioso2018projected}, it is convenient to design the preconditioning matrix as
\begin{equation}\label{eq:phi_matrix}
\Phi = \left[\begin{array}{ccc}
\diag(\{ \alpha_i^{-1} I_{n^2_i} \}_{i \in \mc{I}}) & -A^\top & -D^\top\\
-A & \gamma^{-1} I_m & 0\\
-D & 0 & \tau^{-1} I_p
\end{array}
\right],
\end{equation}
where every $\alpha_i$, $i \in \mc{I}$, $\gamma$ and $\tau$ are chosen large enough to guarantee $\Phi \succ 0$. Essentially, the off-diagonal terms of $\Phi$ in \eqref{eq:phi_matrix} allow to cancel the skew-symmetric matrices in \eqref{eq:low_term}, while the step sizes constitute the main diagonal.

\begin{algorithm}[!t]
	\caption{Semi-decentralized \gls{pFB}}\label{alg:pFB_semi}
	\DontPrintSemicolon
	\SetArgSty{}
	\SetKwFor{ForAll}{for all}{do}{end forall}
	\textbf{Initialization:}
	\begin{itemize}
		\item[] $\bs{x}^{0} \in \mc{S}$, $\nu^{0} \in \R^{m}$, $\lambda^{0} \in \R^p$
	\end{itemize}
	\smallskip
	\textbf{Iteration $(k \in \N)$:}
	\begin{itemize}
		\item[] $\bs{x}^{k+1} = \proj_{\mc{S}} \left( \bs{x}^{k} - \alpha \left(F + A^\top \nu^{k} + D^\top \lambda^{k} \right) \right)$
		\item[]\vspace{-.25cm}
		\item[] $\nu^{k+1} = \nu^k + \gamma \left(2 A \bs{x}^{k+1} - A \bs{x}^{k} - b \right)$
		\item[]\vspace{-.25cm}
		\item[] $\lambda^{k+1} = \lambda^k + \tau \left(2 D \bs{x}^{k+1} - D \bs{x}^{k}\right)$
	\end{itemize}
\end{algorithm}

The steps in Algorithm~\ref{alg:pFB_semi} summarize the proposed semi-decentralized \gls{pFB}.
Remarkably, such a procedure is suitable for a parallel implementation. Specifically, each agent $i \in \mc{I}$ first updates its decision variables $x_i$ by means of a projection on a local \gls{psd} cone of reduced dimension, i.e.,
$$
\begin{aligned}
x&_i^{k + 1} = \proj_{\mc{S}_i} \left( x_i^k - \alpha_i \left( c_i + A_i^\top \nu^k + D_i^\top \lambda^k \right) \right)\\ 
&= \vvec\left( \proj_{\bbS^{n_i}_{\succcurlyeq 0}} \left( \mat\left( x_i^k - \alpha_i \left( c_i + A_i^\top \nu^k + D_i^\top \lambda^k \right) \right) \right) \right),
\end{aligned}
$$
where the projection onto $\bbS^{n_i}_{\succcurlyeq 0}$ of $\mat\left( x_i^k - \alpha_i \left( c_i + A_i^\top \nu^k + D_i^\top \lambda^k \right) \right)$ can be efficiently computed via eigenvalue decomposition. Then, the agents communicate their decisions to a central coordinator, which is in charge of updating both the Lagrange multipliers associated with the consistency constraints ($\lambda$), and the dual variable associated with the original coupling constraints ($\nu$). Note that the central coordinator performs almost inexpensive computations.
Due to privacy reasons, we next investigate how to avoid this direct involvement by a central coordinator, by proposing a fully distributed \gls{pFB} algorithm to solve the \glspl{SDP} in \eqref{eq:SDP_single_prob_vecform}.

\subsection{Distributed forward-backward algorithm}
As in \cite{yi2019operator}, to design a distributed algorithm based on the \gls{pFB} splitting, we endow each agent with a local copy of dual variables $\nu_{i} \in \R^{m}$ and $\lambda_{i} \in \R^{p}$. Moreover, we augment the state of each agent $i \in \mc{I}$ with two local auxiliary variables, namely $z_i \in \R^{m}$ and $y_i \in \R^{p}$, with the aim to enforce local consensus on the copies of the multipliers. 

Let us now preliminary introduce some useful quantities, such as $L_{\nu} \coloneqq L \otimes I_m$, $L_{\lambda} \coloneqq L \otimes I_p$, $\Lambda \coloneqq \diag(\left\{A_i\right\}_{i \in \mc{I}})$, $\bar{b} \coloneqq b\otimes\bsone_{N}$, and $\Xi \coloneqq \diag(\left\{D_i\right\}_{i \in \mc{I}})$. Here, $L$ is the weighted Laplacian matrix associated to $\mc{G}$. Furthermore, $\bar{z} \coloneqq \col(\left\{z_i\right\}_{i \in \mc{I}}) \in \R^{mN}$, $\bar{y} \coloneqq \col(\left\{y_i\right\}_{i \in \mc{I}}) \in \R^{pN}$, while $\bar{\nu}$ and $\bar{\lambda}$ are defined similarly. Then, by augmenting the state and by adding the consensus constraints on the dual variables, the operators in \eqref{eq:A}--\eqref{eq:B} extend as follows:
\begin{align}
\bar{\mc{A}} : &\left[\begin{array}{c}
\bs{x}\\
\bar{z}\\
\bar{y}\\
\bar{\nu}\\
\bar{\lambda}
\end{array}\right] \mapsto  \left[\begin{array}{c}
F\\
0\\
0\\
\bar{b} + L_{\nu} \bar{\nu}\\
L_{\lambda} \bar{\lambda}
\end{array}\right],\label{eq:A_bar}\\
\bar{\mc{B}} : &\left[\begin{array}{c}
\bs{x}\\
\bar{z}\\
\bar{y}\\
\bar{\nu}\\
\bar{\lambda}
\end{array}\right] \mapsto \left[\begin{array}{ccccc}
\phantom{-}0 & 0 & 0 & \phantom{-}\Lambda^\top & \phantom{-}\Xi^\top\\
\phantom{-}0 & 0 & 0 & -L_{\nu} & \phantom{-}0\\
\phantom{-}0 & 0 & 0 & \phantom{-}0 & -L_{\lambda}\\
- \Lambda & L_{\nu} & 0 & \phantom{-}0 & \phantom{-}0\\
- \Xi & 0 & L_{\lambda} & \phantom{-}0 & \phantom{-}0\\
\end{array}
\right] \left[\begin{array}{c}
\bs{x}\\
\bar{z}\\
\bar{y}\\
\bar{\nu}\\
\bar{\lambda}
\end{array}\right]\nonumber\\
&\hspace{1.2cm}+ \left[\begin{array}{c}
\textsc{N}_{\mc{S}}(\bs{x})\\
0\\
0\\
0\\
0
\end{array}\right].\label{eq:B_bar}
\end{align}

Thus, by defining $\bs{\varpi} \coloneqq \col(\bs{x}, \bar{z}, \bar{y}, \bar{\nu}, \bar{\lambda})$, we mimic the steps for the semi-decentralized case to design a suitable preconditioning matrix $\bar{\Phi}$ so that the mapping $\left(\textrm{Id} + \bar{\Phi}^{-1} \bar{\mc{B}}\right)^{-1} \circ \left(\textrm{Id} - \bar{\Phi}^{-1} \bar{\mc{A}}\right)$ can be iterated. Also in this case, we have that $\bs{\varpi} = (\textrm{Id} + \bar{\Phi}^{-1} \mc{B})^{-1} \circ \left(\textrm{Id} - \bar{\Phi}^{-1} \mc{A}\right) \bs{\varpi}$ if and only if $\bs{\varpi} \in \zer(\bar{\mc{A}} + \bar{\mc{B}}) = \zer(\bar{\Phi}^{-1} (\bar{\mc{A}} + \bar{\mc{B}}))$. Specifically, by resorting on the inclusion $- \bar{\mc{A}} \left( \bs{\varpi}^{k} \right) \in \bar{\Phi} \left( \bs{\varpi}^{k+1} - \bs{\varpi}^{k}\right) + \bar{\mc{B}} \left( \bs{\varpi}^{k+1} \right)$, after substituting the operators $\bar{\mc{A}}$ and $\bar{\mc{B}}$ in \eqref{eq:A_bar}--\eqref{eq:B_bar}, the matrix $\bar{\Phi}$ can be conveniently chosen as
\begin{equation}\label{eq:phi_bar_matrix}
\bar{\Phi} = \left[\begin{array}{ccccc}
\phantom{-}\alpha^{-1} & 0 & 0 & -\Lambda^\top & -\Xi^\top\\
\phantom{-}0 & \bar{\sigma}^{-1} & 0 & L_{\nu} & 0\\
\phantom{-}0 & 0 & \bar{\eta}^{-1} & 0 & L_{\lambda}\\
-\Lambda & L_{\nu} & 0 & \bar{\gamma}^{-1} & 0\\
-\Xi & 0 & L_{\lambda} & 0 & \bar{\tau}^{-1}\\
\end{array}
\right].
\end{equation}

Also in this case, the entries on the main diagonal, with $\bar{\sigma} \coloneqq \diag(\{\sigma_i\}_{i \in \mc{I}}) \otimes I_m$ ($\bar{\eta}$, $\bar{\gamma}$ and $\bar{\tau}$ are defined similarly), are chosen large enough to have a positive definite $\bar{\Phi}$. Here, the off-diagonal terms of $\bar{\Phi}$ in \eqref{eq:phi_bar_matrix} cancel the skew-symmetric terms in \eqref{eq:B_bar} and the step sizes are located on the main diagonal. The distributed \gls{pFB} procedure is summarized in Algorithm~\ref{alg:pFB_distr}, where the steps to be performed at every iteration by each single agent are emphasized. Specifically, the algorithm alternates communication and computation tasks. Once that the agent receives the local copy of the dual variables from its neighbors in $\mc{N}_i$, it updates the augmented state, i.e., its decision variable $x_i$ together with $z_i$ and $y_i$ (\texttt{S1}). Then, the agents transmit to their neighbors the updated version of the augmented state to compute a local update of the copies of the multipliers (\texttt{S2}). Finally, we emphasize that each computation is based on either local or communicated quantities only, and it does not require a central coordinator.
\begin{algorithm}[!t]
	\caption{Distributed \gls{pFB}}\label{alg:pFB_distr}
	\DontPrintSemicolon
	\SetArgSty{}
	\SetKwFor{ForAll}{for all}{do}{end forall}
	\textbf{Initialization:} For each agent $i \in \mc{I}$
	\begin{itemize}
		\item[] $x_i^{0} \in \mc{S}_i$, $z_i^0, \, \nu_i^{0} \in \R^{m}$, $y^0_i, \, \lambda_i^{0} \in \R^p$	
	\end{itemize}
	\smallskip
	\textbf{Iteration $(k \in \N)$:} For each agent $i \in \mc{I}$\\
	\texttt{S1)} Receives $\nu^k_{j}$ and $\lambda^k_{j}$, $j \in \mc{N}_i$, and updates
	\begin{itemize}
		\item[] $x_i^{k+1} = \proj_{\mc{S}_i} \left( x_i^{k} - \alpha_i \left(c_i + A_i^\top \nu_i^{k} + D_i^\top \lambda_i^{k} \right) \right)$
		\item[]\vspace{-.25cm}
		\item[] $z_i^{k+1} = z_i^k + \sigma_i \sum_{j \in \mc{N}_i} w_{i,j} (\nu^k_{i} - \nu^k_{j})$
		\item[]\vspace{-.25cm}
		\item[] $y_i^{k+1} = y_i^k + \eta_i \sum_{j \in \mc{N}_i} w_{i,j} (\lambda^k_{i} - \lambda^k_{j})$
	\end{itemize}
	\smallskip
	\texttt{S2)} Receives $z^{k+1}_{j}$ and $y^{k+1}_{j}$, $j \in \mc{N}_i$, and updates
	\begin{itemize}
		\item[] $\nu_i^{k+1} = \nu_i^k + \gamma_i [ A_i (2 x_i^{k+1} - x_i^{k}) - b_i$\\ 
		\hspace{1.1cm}$- \sum_{j \in \mc{N}_i} w_{i,j} ( 2 ( z_i^{k+1} - z_j^{k+1}) - (z_i^{k} - z_j^{k}))$\\ 
		\hspace{1.1cm}$- \sum_{j \in \mc{N}_i} w_{i,j} (\nu_i^{k} - \nu_j^{k}) ]$
		\item[]\vspace{-.25cm}
		\item[] $\lambda_i^{k+1} = \lambda_i^k + \tau_i [ D_i (2 x_i^{k+1} - x_i^{k})$\\
		\hspace{1.1cm}$- \sum_{j \in \mc{N}_i} w_{i,j} ( 2 ( y_i^{k+1} - y_j^{k+1}) - (y_i^{k} - y_j^{k}))$\\ 
		\hspace{1.1cm}$- \sum_{j \in \mc{N}_i} w_{i,j} (\lambda_i^{k} - \lambda_j^{k}) ]$
	\end{itemize}
\end{algorithm}

\section{Numerical simulations}

\begin{figure}
	\centering
	\includegraphics[width=.8\columnwidth]{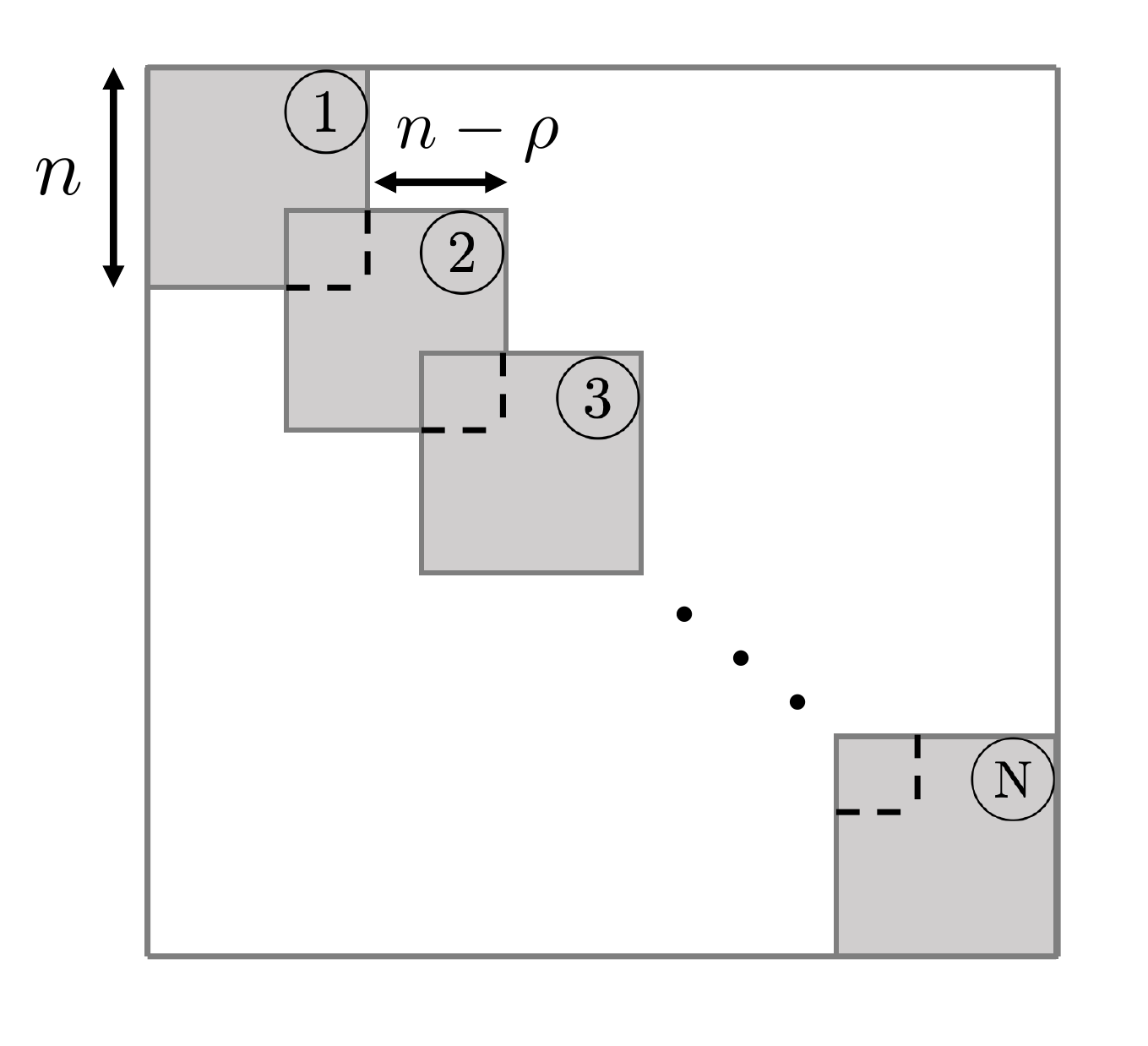}
	\caption{Banded sparsity pattern.}
	\label{fig:banded_sparsity}
\end{figure}

We first compare the performances of the presented algorithms on randomly generated \gls{SDP} problems. The simulations are run in Matlab environment on a laptop with an Intel Core i5 2.9 GHz CPU and 8 Gb RAM.
Specifically, we focus on an \gls{SDP} characterized by a banded sparsity pattern illustrated in Fig.~\ref{fig:banded_sparsity}. Here, each agent $i \in \mc{I} \coloneqq \{1,\ldots,N\}$ controls a diagonal block of dimension $n_i \times n_i$, while sharing a subset of variables with the $(i+1)$-th agent (the parameter $\rho > 0$ specifies the number of overlapping elements). For simplicity, we assume $n_i = n$ for all $i \in \mc{I}$, and therefore the associated \gls{psd} cone has dimension $nN - \rho(N-1)$. As in \cite{zheng2019chordal}, we first generate random symmetric matrices $A_k$, $k \in \mc{K}$, and $W$ with banded sparsity and entries drawn uniformly from $(0,1)$. Then, we generate a strictly feasible matrix $X_\textrm{f} = W + \kappa I$, with $\kappa > 0$ large enough to guarantee $X_\textrm{f} \succ 0$. Finally, the vector $b$ of equality constraints is computed as $b_k = \inner{A_k}{X_\textrm{f}}$, $k \in \mc{K}$.
\begin{figure}
	\centering
		\includegraphics[width=.9\columnwidth]{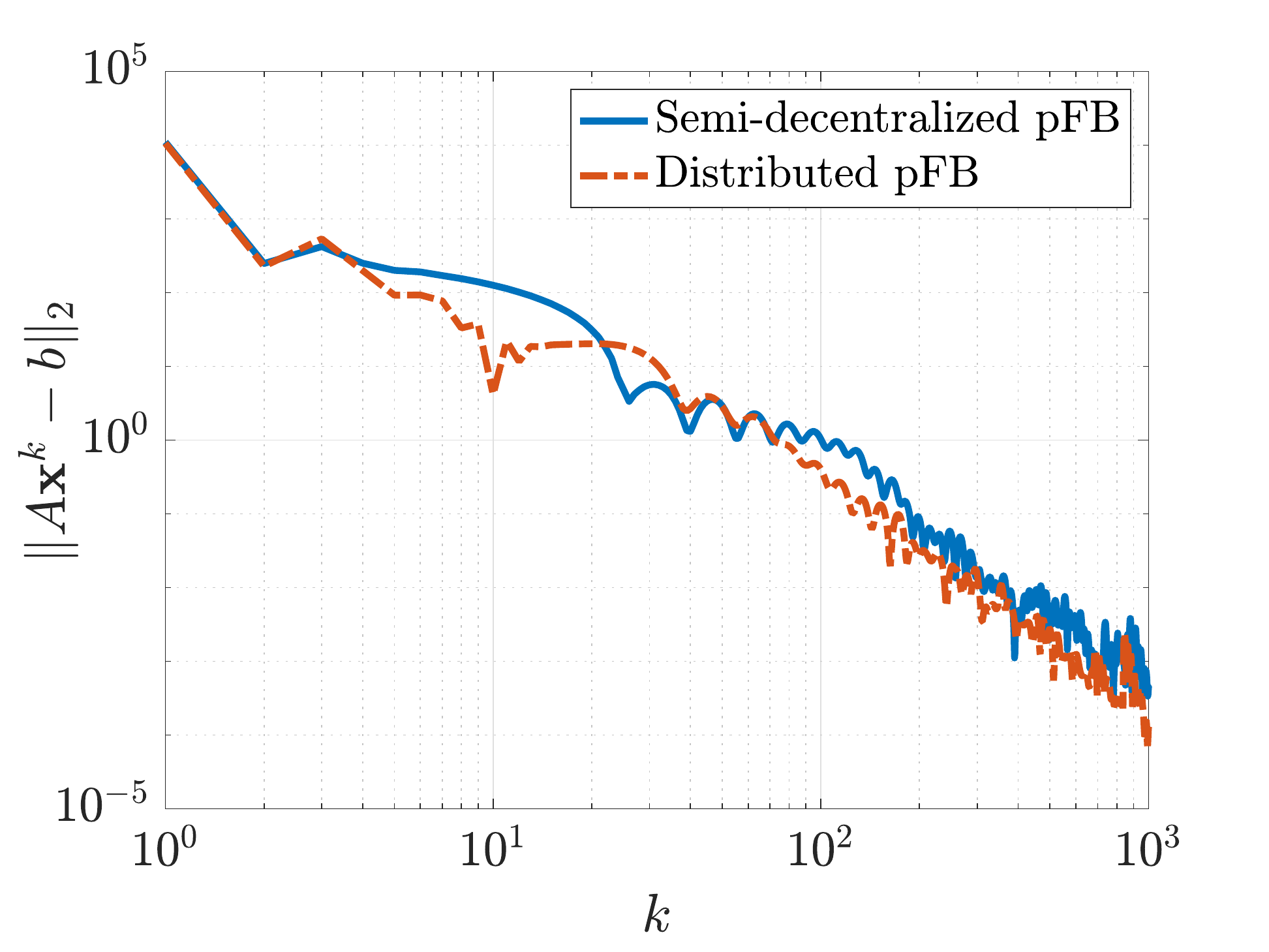}
	\caption{Coupling equality constraints.}
	\label{fig:norm_cons}
\end{figure}
\begin{figure}
	\centering
		\includegraphics[width=.9\columnwidth]{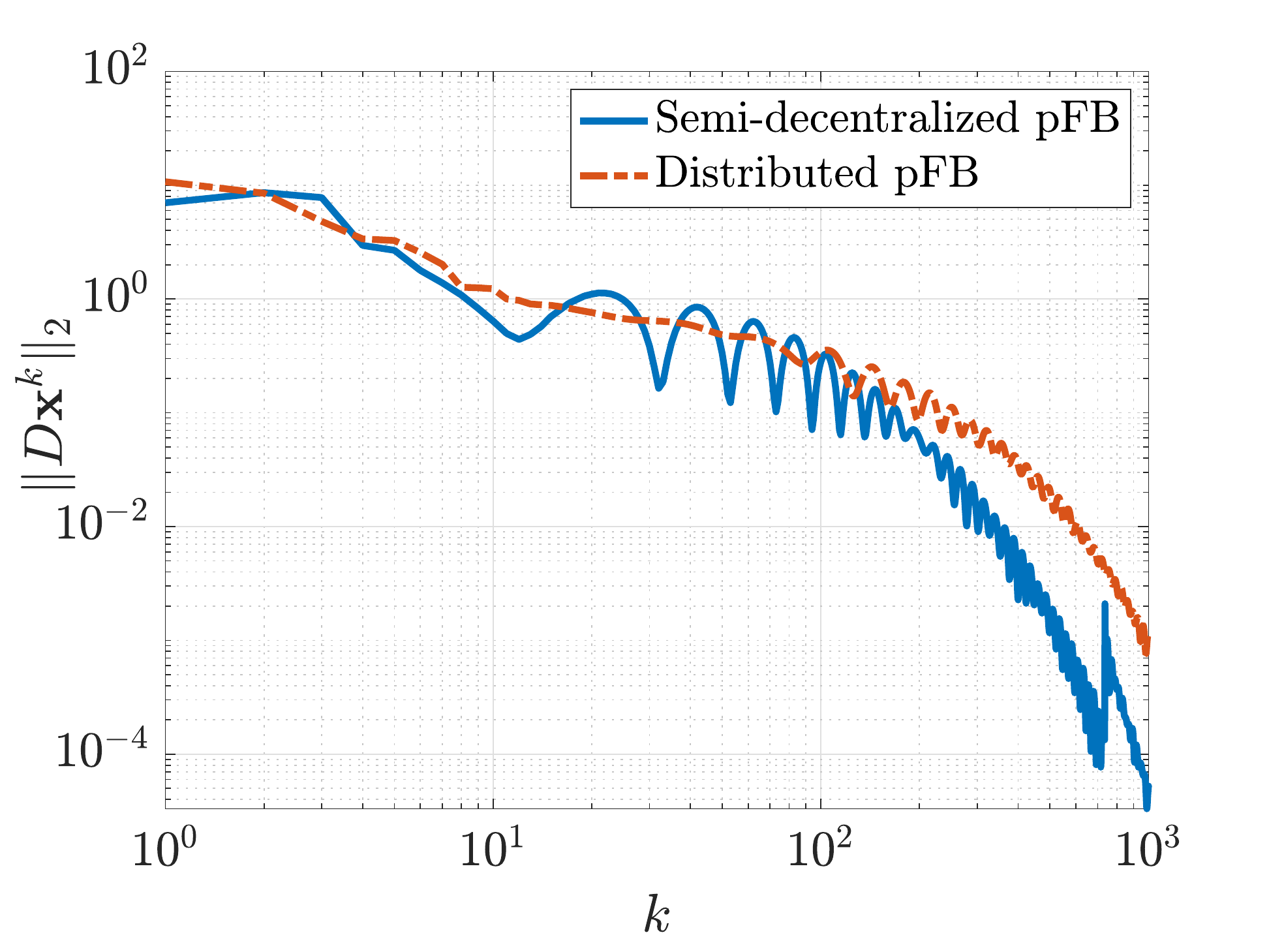}
	\caption{Consistency constraints.}
	\label{fig:cons_cons}
\end{figure}

\begin{figure*}
	\centering
	\includegraphics[width=\textwidth]{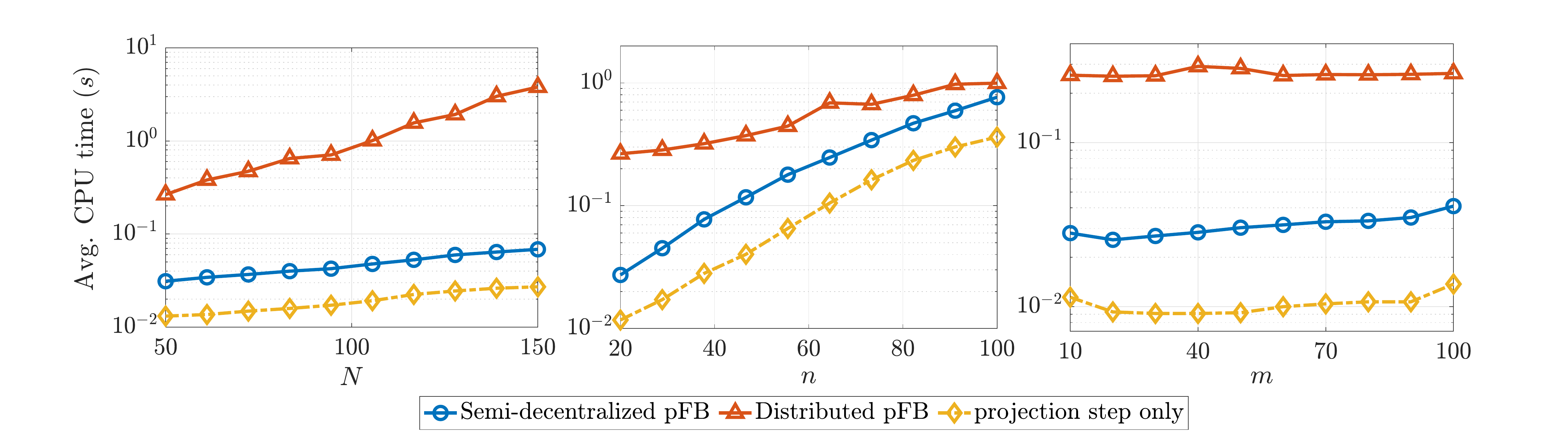}
	\caption{Average CPU time per $100$ iterations. From left to right: varying the number of diagonal blocks (i.e., agents); varying the dimension of each block; varying the number of coupling constraints of the original \gls{SDP}.}\label{fig:cpu}
\end{figure*}

\begin{figure}
	\centering
	\includegraphics[width=.9\columnwidth]{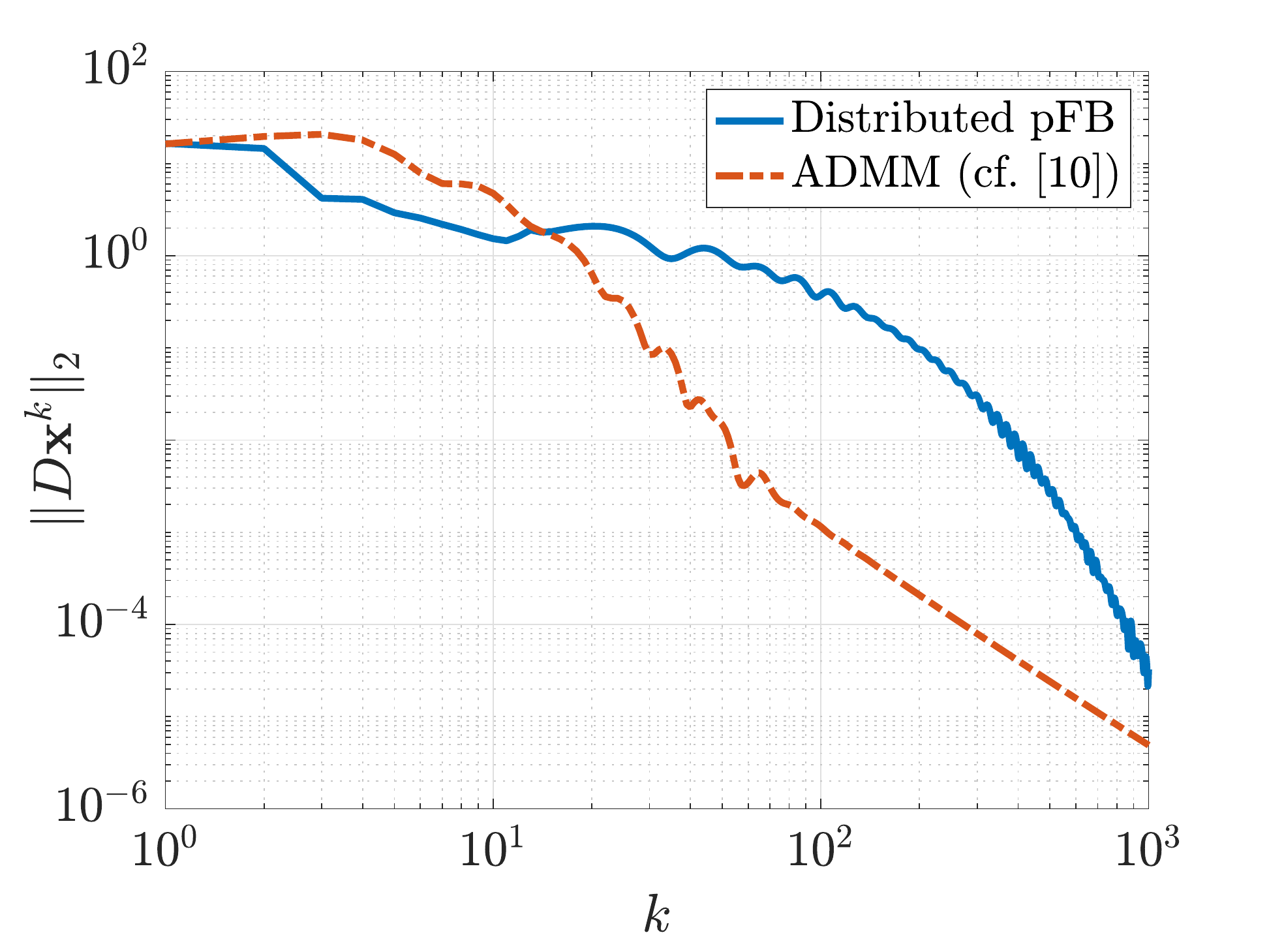}
	\caption{Comparison of the convergence behavior of the consistency constraints between Algorithm~\ref{alg:pFB_distr} and the ADMM algorithm in \cite{kalbat2015fast}. }\label{fig:comparison_consistency_cons}
\end{figure}

\begin{figure}
	\centering
	\includegraphics[width=.9\columnwidth]{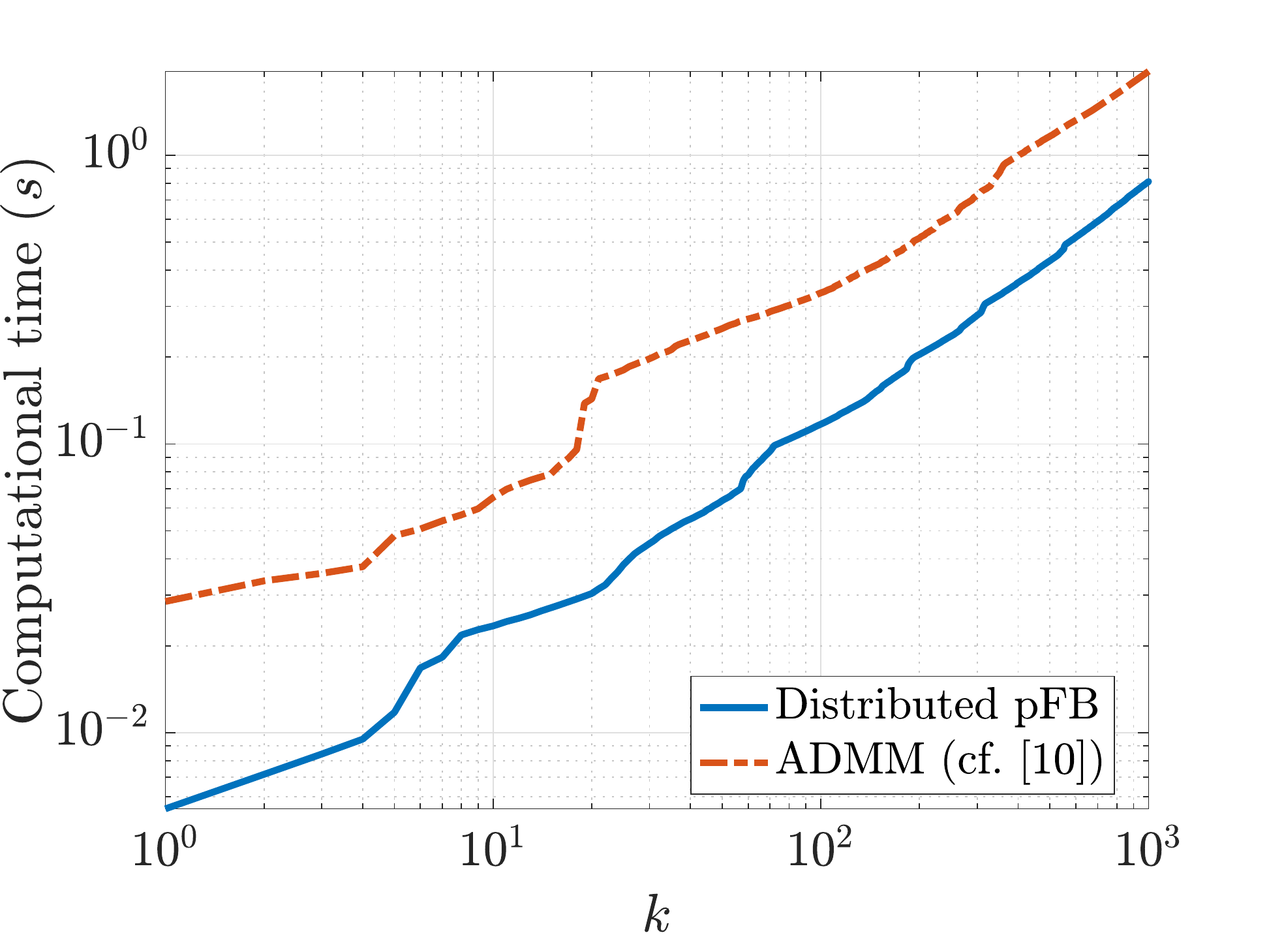}
	\caption{Comparison of the integral value of the computational time between Algorithm~\ref{alg:pFB_distr} and the ADMM algorithm in \cite{kalbat2015fast}.}\label{fig:IoCPT}
\end{figure}

\begin{figure}
	\centering
	\includegraphics[width=.9\columnwidth]{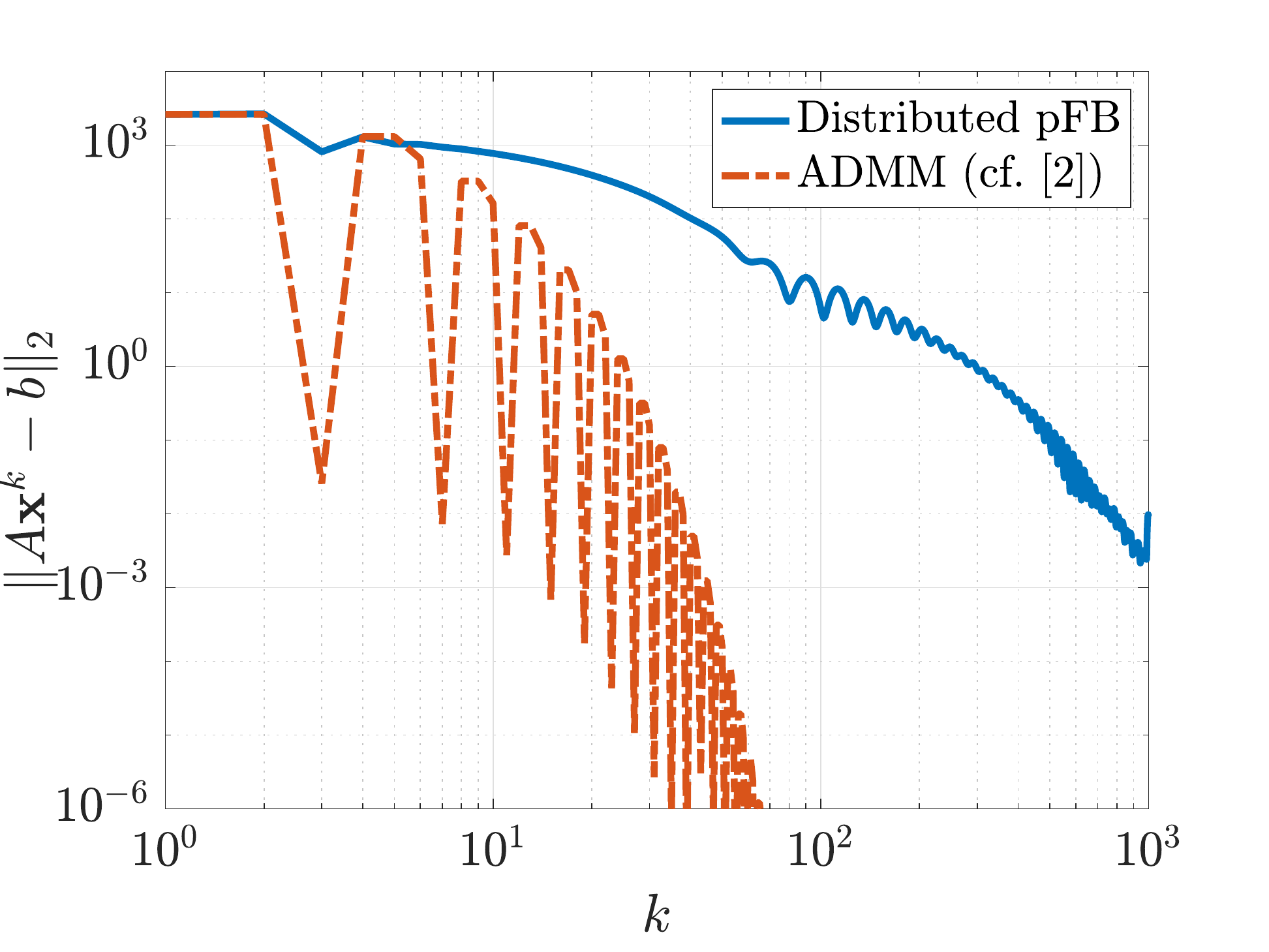}
	\caption{Comparison of the convergence behavior of the equality constraints between Algorithm~\ref{alg:pFB_distr} and the ADMM algorithm in \cite{simonetto2014distributed}.}\label{fig:comparison_SL_equality_norms}
\end{figure}

\begin{figure}
	\centering
	\includegraphics[width=.9\columnwidth]{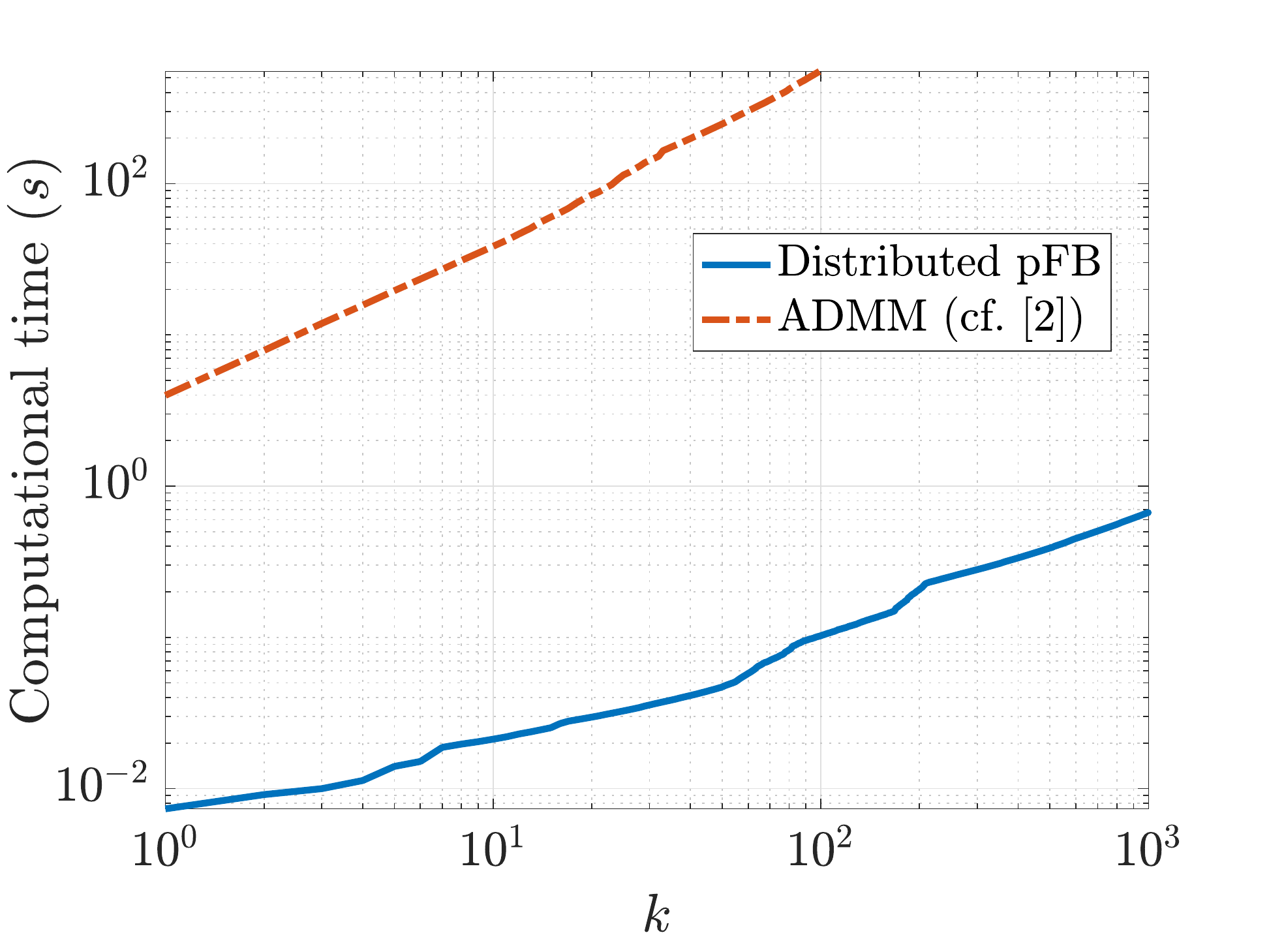}
	\caption{Comparison of the integral value of the computational time between Algorithm~\ref{alg:pFB_distr} and the ADMM algorithm in \cite{simonetto2014distributed}.}\label{fig:comparison_SL_IoCPT}
\end{figure}

By considering an instance with $N = 50$, $n = 20$, $\rho = m = 5$, Figures~\ref{fig:norm_cons}--\ref{fig:cons_cons} show the convergence of both coupling and consistency constraints when applying Algorithm~\ref{alg:pFB_semi} and \ref{alg:pFB_distr}, while Fig.~\ref{fig:cpu} illustrates the average CPU time per $100$ iterations as a function of $N$, $n$ and $m$, respectively. The (averaged) CPU time to compute the projection onto $\bbS^{n_i}_{\succcurlyeq 0}$, for all $i \in \mc{I}$ is also reported. For each plot, only one parameter at time changes, while the other are fixed to $N = 50$, $n = 20$, $m = 10$ and $\rho = 10$.
While the two algorithms exhibit the same behavior in term of convergence, Algorithm~\ref{alg:pFB_semi} takes less time on average to perform each iteration.

Finally, Figures.~\ref{fig:comparison_consistency_cons}--\ref{fig:comparison_SL_IoCPT} compare the behavior of the distributed \gls{pFB} in Algorithm~\ref{alg:pFB_distr} and the ADMM algorithms proposed in \cite{kalbat2015fast} and \cite{simonetto2014distributed} on two different instances with $N = 50$, $n = 30$, $m = 0$ and $\rho = 10$, and $N = 20$, $n = 20$, $m = 10$ and $\rho = 5$, respectively. The ADMM algorithms converge with fewer iterations (Fig.~\ref{fig:comparison_consistency_cons} and \ref{fig:comparison_SL_equality_norms}), at the price of much higher computational cost at each iteration, as shown in Fig.~\ref{fig:IoCPT} and \ref{fig:comparison_SL_IoCPT}. Specifically, since the iterations of the considered algorithms shall be performed in parallel, let $t^k_i$ be the computational time taken by the $i$-th agent, $i \in \mc{I}$, to perform the steps of a certain algorithm at iteration $k \in \N$. The value at the generic iteration $K$ in Fig.~\ref{fig:IoCPT} and \ref{fig:comparison_SL_IoCPT} is computed as $\sum_{k = 1}^{K} \textrm{max}_{i \in \mc{I}} \; t_i^k$. However, we emphasize that the distributed ADMM in \cite{kalbat2015fast} is less general than the \gls{pFB} in Algorithm~\ref{alg:pFB_distr}, since it is tailored for separable \glspl{SDP} with local equality constraints only (this justifies $m = 0$ for the first numerical simulation), while the one in \cite{simonetto2014distributed}, at every iteration $k \in \N$, requires to solve an \gls{SDP} for each $i \in \mc{I}$.

\section{Conclusion and outlook}
Monotone operator splitting theory promises to be a general purpose tool to design scalable and distributed algorithms for decomposable, large-scale \glspl{SDP}. The operator-theoretic perspective given in this paper paves the way for a systematic development of efficient methods with reasonable computational burden and convergence guarantees. Accelerated strategies, as well as additional operator splitting methods, are left as future work.



\bibliographystyle{IEEEtran}
\bibliography{20_ECC.bib}

\end{document}